\documentclass[11pt]{article}

\usepackage{amscd, amsmath, amssymb, fancyhdr, graphicx, epsfig, url, mathrsfs}
\usepackage{graphicx}
\usepackage{epigraph}
\usepackage[backref=page]{hyperref}
\renewcommand*{\backrefalt}[4]{%
	\ifcase #1 (Not cited.)%
	\or        (Cited on page~#2.)%
	\else      (Cited on pages~#2.)%
	\fi}

\hypersetup{
	colorlinks   = true,
	citecolor    = magenta,
	linkcolor    = blue,
	urlcolor     = magenta	
}

\numberwithin{equation}{section}

\newcommand{\version}{version 2.0,\ \ April 18, 2025}

\renewcommand{\leftmark}{{\scriptsize K\"ahler structures on degenerate twistors}}

\makeatletter
\def\x@arrow{\DOTSB\Relbar}
\def\xlongrightarrowfill@{\arrowfill@\relbar\relbar\longrightarrow}
\newcommand{\xlongrightarrow}[2][]{%
        \ext@arrow 0099\xlongrightarrowfill@{#1}{#2}}
\makeatother

\def\eqref#1{(\ref{#1})}

\newcommand{\arrow}{{\:\longrightarrow\:}}

\def\C{{\Bbb C}}

\newcommand{\R}{{\Bbb R}}

\newcommand{\6}{\partial}
\def\1{\sqrt{-1}\:}
\newcommand{\restrict}[1]{|_{#1}}
\newcommand{\cntrct}                
{\hspace{2pt}\raisebox{1pt}{\text{$\lrcorner$}}\hspace{2pt}}

\renewcommand{\tilde}{\widetilde}
\renewcommand{\bar}{\overline}
\renewcommand{\phi}{\varphi}
\renewcommand{\epsilon}{\varepsilon}
\renewcommand{\geq}{\geqslant}

\newcommand{\Teich}{\operatorname{\sf Teich}}

\newcommand{\Comp}{\operatorname{\sf Comp}}

\newcommand{\Per}{\operatorname{\sf Per}}

\newcommand{\Id}{\operatorname{Id}}

\newcommand{\Vol}{\operatorname{Vol}}

\newcommand{\Diff}{\operatorname{\sf Diff}}

\renewcommand{\Re}{\operatorname{Re}}
\renewcommand{\Im}{\operatorname{Im}}

\newcommand{\bbA}{\mathbb{A}}
\newcommand{\bbZ}{\mathbb{Z}}
\newcommand{\bbQ}{\mathbb{Q}}
\newcommand{\bbR}{\mathbb{R}}
\newcommand{\bbC}{\mathbb{C}}
\newcommand{\bbH}{\mathbb{H}}

\newcommand{\bbP}{\mathbb{P}}

\newcommand{\CC}{\mathcal{C}}
\newcommand{\DD}{\mathcal{D}}
\newcommand{\EE}{\mathcal{E}}

\newcommand{\MM}{\mathcal{M}}

\newcommand{\XX}{\mathcal{X}}
\newcommand{\YY}{\mathcal{Y}}

\renewcommand{\ge}{\geqslant}
\renewcommand{\le}{\leqslant}

\newcommand{\st}{\enskip |\enskip}

\newcommand{\cnv}{\,\lrcorner\,}

\newcommand{\dd}{\partial}

\newcommand{\ii}{\sqrt{-1}}

\newcommand{\wdg}{\wedge}

\newcommand{\hrarr}{\hookrightarrow}

\newcommand{\proof}{\noindent{\bf Proof:\ }}

\newcounter{Mycounter}[section]
\newcounter{lemma}[section]
\renewcommand{\thelemma}{{Lemma \thesection.\arabic{lemma}}}
\newcommand{\lemma}{%
    \setcounter{lemma}{\value{Mycounter}}
    \refstepcounter{lemma}
    \stepcounter{Mycounter}
    {\noindent \bf \thelemma:\ }}

\newcounter{claim}[section]
\newcounter{sublemma}[section]
\newcounter{corollary}[section]
\newcounter{theorem}[section]
\renewcommand{\thetheorem}{{Theorem \thesection.\arabic{theorem}}}
\newcommand{\theorem}{%
    \setcounter{theorem}{\value{Mycounter}}
    \refstepcounter{theorem}
    \stepcounter{Mycounter}
    {\noindent \bf \thetheorem:\ }}

\newcounter{conjecture}[section]
\renewcommand{\theconjecture}{{Conjecture \thesection.\arabic{conjecture}}}
\newcommand{\conjecture}{%
    \setcounter{conjecture}{\value{Mycounter}}
    \refstepcounter{conjecture}
    \stepcounter{Mycounter}
    {\noindent \bf \theconjecture:\ }}

\newcounter{proposition}[section]
\renewcommand{\theproposition}
      {{Proposition \thesection.\arabic{proposition}}}
\newcommand{\proposition}{%
    \setcounter{proposition}{\value{Mycounter}}
    \refstepcounter{proposition}
    \stepcounter{Mycounter}
    {\noindent \bf \theproposition:\ }}

\newcounter{definition}[section]
\renewcommand{\thedefinition}
      {{Definition~\thesection.\arabic{definition}}}
\newcommand{\definition}{%
    \setcounter{definition}{\value{Mycounter}}
    \refstepcounter{definition}
    \stepcounter{Mycounter}
    {\noindent \bf \thedefinition:\ }}

\newcounter{example}[section]
\newcounter{remark}[section]
\renewcommand{\theremark}{{Remark \thesection.\arabic{remark}}}
\newcommand{\remark}{%
    \setcounter{remark}{\value{Mycounter}}
    \refstepcounter{remark}
    \stepcounter{Mycounter}
    {\noindent \bf \theremark:\ }}

\newcounter{problem}[section]
\newcounter{question}[section]
\makeatletter

\@addtoreset{equation}{section} \@addtoreset{footnote}{section}
\makeatother

\def\blacksquare{\hbox{\vrule width 5pt height 5pt depth 0pt}}
\def\endproof{\blacksquare}

\begin{document}

\begin{center}
{\LARGE\bf
Hermitian-symplectic and\\
K\"ahler structures on degenerate twistor deformations
}

\bigskip

Andrey Soldatenkov, 
Misha Verbitsky\footnote{Partially supported by 
FAPERJ SEI-260003/000410/2023 and CNPq - Process 310952/2021-2. \\


 }

\end{center}

{\small \hspace{0.15\linewidth}
\begin{minipage}[t]{0.7\linewidth}
{\bf Abstract} 
Let $(M, \Omega)$ be a holomorphically symplectic
manifold equipped with a holomorphic Lagrangian fibration
$\pi: M \to B$, and $\eta$ a closed $(1,1)$-form on $B$.
Then $\Omega+ \pi^* \eta$ is a holomorphically symplectic
form on a complex manifold which is called 
the degenerate twistor deformation of $M$.
We prove that degenerate twistor deformations 
of compact holomorphically symplectic K\"ahler
manifolds are also K\"ahler. First, we prove that
degenerate twistor deformations are Hermitian
symplectic, that is, tamed by a symplectic form;
this is shown using positive currents and 
an argument based on the Hahn--Banach theorem,
originally due to Sullivan. Then we apply a version of Huybrechts's
theorem showing that two non-separated 
points in the Teichm\"uller space of holomorphically
symplectic manifolds correspond to bimeromorphic
manifolds if they are Hermitian symplectic.
\end{minipage}
}

\tableofcontents


\section{Introduction}


\subsection{Degenerate twistor deformations and the Tate--Shafare\-vich twist}

Let $(M, I, \Omega)$ be a holomorphically symplectic
manifold, that is, an almost complex manifold equipped with
a closed, non-degenerate $(2,0)$-form. 
Clearly, the almost complex structure $I$ is uniquely
determined by the differential form $\Omega$. It is not hard
to see that the integrability of the almost complex structure $I$
follows from closedness of $\Omega$ (\cite{_V:degenerate_}).
This allows one to define the complex structure in terms
of $\Omega$, giving a motivation for the notion of
C-symplectic form (\ref{_C_symp_Definition_}).

The notion of degenerate twistor deformation,
introduced in \cite{_V:degenerate_},
and further explored in \cite{BDV} and \cite{SV},
is one of the applications of this approach.
Given a holomorphic Lagrangian fibration $\pi:\; M \arrow B$
on a holomorphic symplectic manifold $(M, \Omega)$, we can 
replace $\Omega$ by $\Omega+t\pi^*\eta$, where $\eta$ 
is a closed $(2,0)+(1,1)$-form on $B$, and $t\in \C$.
The new form $\Omega_t:=\Omega+t\pi^*\eta$ is C-symplectic,
that is, holomorphic symplectic with respect to
a complex structure $I_t$ (\ref{thm_csympl_lagr}).
This defines a holomorphic family of complex manifolds $(M,I_t), t\in \C$,
which is called {\bf the degenerate twistor deformation}.

It shares many of its properties with twistor deformations
(the ones which come from a hyperk\"ahler rotation);
indeed, the degenerate twistor deformation can be obtained
as a limit of twistor deformations (\cite[Section 1.2]{_V:degenerate_}).

A version of this construction is known as the ``Tate--Shafarevich twist''.
For elliptic surfaces, it was explored at some length in 
\cite{_FM:4-manifolds_}, and for higher dimensions in 
\cite{_Markman:Lagrangian_} and \cite{_AbashevaRogov_}.
In \cite{_AbashevaRogov_} it was shown that
the Tate--Shafarevich twist and the degenerate twistor deformation
coincide if the Lagrangian fibration has no multiple fibers
in codimension 1. 

In the present paper, we prove that 
the degenerate twistor deformation is always K\"ahler (\ref{thm_main}), 
answering the question which was originally posed by E. Markman.

\hfill

\theorem
Let $M$ be a compact holomorphically symplectic
manifold of K\"ahler type, $\pi:\; M \arrow B$ a Lagrangian
fibration, and $\tau:\; {\cal M}\arrow \C$ the corresponding degenerate twistor deformation.
Then all fibers of $\tau$ admit a K\"ahler metric.

\proof
\ref{thm_main}.
\endproof

\hfill

In \cite{_Markman:Lagrangian_}, Markman has studied
the Tate--Shafarevich twist of a Lagrangian fibration on a
Hilbert scheme of a K3 surface. Markman calls a compact
holomorphically symplectic manifold $(M, \Omega)$ 
``M-special'' if the space $\langle [\Re\Omega],[\Im\Omega]\rangle \subset H^2(M, \R)$
does not contain non-zero rational cohomology classes. He has shown that
the Tate--Shafarevich twist of a Lagrangian fibration on a
Hilbert scheme of a K3 surface is K\"ahler, if it is not
M-special. In \cite{_AbashevaRogov_}, this result was generalized to all Lagrangian
fibrations without multiple fibers in codimension 1.
In the present paper we give an independent proof that works for
arbitrary compact hyperk\"ahler manifolds without the
assumption of non-M-specialty.

\subsection{Non-K\"ahler deformations of hyperk\"ahler manifolds}

The earliest publication on K\"ahlerness of holomorphically
symplectic manifolds was a preprint of Andrey Todorov
\cite{_Todorov:MPIM_}. Using a Harvey-Lawson-type Hahn-Banach argument,
Todorov tried to show that a simply connected, compact,
holomorphically symplectic manifold with $H^{2,0}(M)=\C$
is K\"ahler.

The earliest counterexamples to this statement were obtained in 1994
by D. Guan, \cite{_Guan:1_,_Guan:2_,_Guan:3_,_Bogomolov:Guan_}.
These are manifolds, obtained from the Kodaira surface
(non-Kahler, holomorphically symplectic complex surface)
by a procedure reminding of the one which is used to construct
the generalized Kummer variety from a 2-dimensional torus.
This manifold is non-K\"ahler because it contains a
complex surface isomorphic to a blow-up of Kodaira surface
(\cite[Corollary 4.10]{_Bogomolov:Guan_}).

Another example of a non-K\"ahler holomorphically
symplectic manifold comes from birational
geometry.\footnote{We benefited from a Mathoverflow
discussion at
\url{http://mathoverflow.net/questions/302649/are-there-non-projective-but-algebraic-hyperkahler-varieties}
and the information provided by the user YangMills.}
Recall that {\bf Mukai flop} is a bimeromorphic transform
which takes a Lagrangian $\C P^n$ in a holomorphically
symplectic manifold $M$, and replaces it with another one,
obtaining another holomorphically symplectic manifold
(say, $M'$). The second homology groups of these manifolds
are isomorphic, however, the homology 
class of a rational curve (a line contained in the $\C P^n$) on $M$ is opposite to that on $M'$.
If we start from a holomorphically symplectic manifold containing two disjoint 
Lagrangian planes with homologous lines, 
and do a flop in only one of these, we will obtain
a holomorphically symplectic manifold which contains
a complex curve (the union of two lines) homologous to zero.
This manifold is Fujiki class C (that is, bimeromorphic to K\"ahler).
For explicit examples of non-K\"ahler Fujiki class C holomorphic
symplectic manifolds and a reference, see \cite[Example 21.9]{_GHJ_}.

In the sequel, {\bf hyperk\"ahler} means ``compact
holomorphically symplectic manifold of K\"ahler type''.
In \cite{Hu}, D. Huybrechts proved that two bimeromorphic
hyperk\"ahler manifolds corresponds to non-separable (in
the sense of being non-Hausdorff) points of the
corresponding Teichm\"uller space. Conversely,
given two non-separable points $I_1, I_2$ in the 
Teichm\"uller space, the hyperk\"ahler manifolds
$(M, I_1)$ and $(M, I_2)$ are bimeromorphic. Curiously
enough, this argument uses the K\"ahlerness in an essential
way.

If one drops the K\"ahlerness assumption (or a weaker
assumption of existence of the Hodge decomposition on $H^2(M)$ as in
\cite{_Perego:Kahlerness_}; see the discussion of this
work below), the proof is still unknown, though we may
dare to make a conjecture.

\hfill

\conjecture\label{_central_f_hk_Conjecture_}
Let $(M_t, \Omega_t)$ be a smooth family of compact holomorphically
symplectic manifolds over a disk $\Delta$. Assume that all fibers
except the central one are hyperk\"ahler. Then the central
fiber $(M_0, \Omega_0)$ is of Fujiki class C.

\hfill

In \ref{_two_fibers_Herm_symple_bimero_Proposition_}, we prove this conjecture when 
the central fiber $M_0$ is Hermitian symplectic.
In \cite{_Perego:Kahlerness_}, A. Perego discussed
\ref{_central_f_hk_Conjecture_} when $M_0$
satisfies $b_2(M_0)=h^{2,0}(M_0)+h^{1,1}(M_0)+h^{0,2}(M_0)$.

Let us indicate where Huybrechts's proof uses the
K\"ahlerness.  Consider two deformations $M_t$ and $M_t'$ of a
hyperk\"ahler manifold, with $t\in \Delta$. Assume that 
there is an isomorphism $M_t\cong M'_t$ for all $t\neq 0$.
The graph of this isomorphism defines a family of
subvarieties $\Gamma_t$ in $M_t \times M'_t, t\neq 0$.
Using Bishop's compactness theorem, Huybrechts proves
that the family $\Gamma_t$ can be extended to $t=0$,
and the corresponding cycle $\Gamma_0$ contains a component that is the graph of
a bimeromorphic map $M_0 \dashrightarrow M'_0$.
When $M_0$ is non-K\"ahler, Bishop's compactness fails
and the argument becomes invalid. However, 
Bishop's compactness still holds for Hermitian symplectic
manifolds (\ref{_Bishop_Herm_Sympl_Proposition_}), and this leads to 
\ref{_two_fibers_Herm_symple_bimero_Proposition_}.

\hfill

The application of this argument often  brings an error, which
was observed by Y.-T. Siu \cite{_Siu:K3_,_Siu:Todorov_MR_} 
when he was discussing a paper of A. Todorov \cite{_Todorov:K3_}
``Applications of the K\"ahler-Einstein-Calabi-Yau metric to 
moduli of K3 surfaces.'' The same problem occurs later 
in two important and very interesting papers \cite{_Perego:Kahlerness_,_Abasheva:2_}.

As explained by Siu, the gaps occur in the application of
E. Bishop's theorem on limits of complex analytic subvarieties.
Suppose that we have two families of holomorphically symplectic manifolds 
over a disk, isomorphic except in the central fiber; as in Huybrechts's proof,
one considers the limit of the graphs $\Gamma_t$ of the isomorphisms between the 
non-central fibers. The idea is to apply Bishop's theorem
to show that this limit is a complex analytic subvariety
with the same fundamental class.

To apply Bishop's theorem, bounding the volume of the
limiting cycle, one needs to choose a Hermitian form 
on the central fiber, such that the volume of $\Gamma_t$ 
with respect to this form stays bounded for all $t$.
Even if we construct a sequence of K\"ahler metrics on the neighboring fibers
whose cohomology classes converge to a cohomology
class on the central fiber (such as in \cite[Lemma 2.5]{_Perego:Kahlerness_}
and in \cite[Lemma 4.6.1]{_Abasheva:2_}), this does not imply
the volume estimates. Indeed, a cohomology class with positive BBF square
on a holomorphic symplectic manifold that is a limit of K\"ahler
classes on neighboring fibers needs not be K\"ahler itself.

\hfill

\ref{_central_f_hk_Conjecture_} 
is a special case of a more ambitious conjecture proposed
in \cite{LiLiu}.

\hfill

\conjecture\label{_central_f_Kahler_Conjecture_}
Let $M_t$ be a smooth family of compact complex
manifolds over a disk $\Delta$. Assume that all fibers
except the central one are K\"ahler. Then the central
fiber $M_0$ is Fujiki class C.

\hfill

In \cite{LiLiu}, this conjecture is proven
under two extra assumptions on a central fiber:
(a) $M_0$ admits a Hermitian form $\omega$
satisfying $dd^c\omega=0$ and 
$\6 \omega \wedge \bar\6\omega=0$, 
and (b) $h^{0,2}(M_0)= h^{0,2}(M_t)$ for
$t\neq 0$.

\hfill

In \cite{_Popovici:limits_},  this conjecture is proven
when $M_t$, $t\neq 0$ are projective, and
$h^{0,1}(M_0)= h^{0,1}(M_t)$ for
$t\neq 0$; in this case, the central fiber is Moishezon.
Later, this result was amplified by
\cite{_Li_Rao_Wang_Wang_}, who proved that a limit
of Moishezon manifolds is Moishezon.

\subsection{Hermitian symplectic geometry and Streets--Tian conjecture}

The original impulse to this paper is
an observation due to Peternell:
any non-projective Moishezon manifold
admits a positive, exact $(n-1, n-1)$-current,
see \cite{_Peternell:algebraicity_}.
We were trying to use this observation
to prove that a degenerate twistor
deformation of a hyperk\"ahler manifold cannot
be non-projective Moishezon, in hope
to extend it further to Fujiki class C,
and then try to prove \ref{_central_f_hk_Conjecture_}
in this particular situation.

In \ref{thm_Hermitian} we accomplish the first task,
indeed, showing that a  degenerate twistor
deformation of a hyperk\"ahler manifold 
cannot support a  positive, exact $(n-1, n-1)$-current.

Let us illustrate this in a simpler situation,
showing that a degenerate twistor
deformation cannot admit a complex curve
$C$ which is homologous to zero.

Let $\pi:\; M \arrow B$ be a Lagrangian fibration
on a hyperk\"ahler manifold $(M, \Omega)$.
Then $B$ is projective; let $\omega\in H^2(B, \R)$
be its K\"ahler class. Since $\int_C \pi^* \omega=0$,
$C$ is contained in the fibers of $\pi$. However,
the degenerate twistor deformation does not
change the complex structures on the fibers, hence
$C$ is a complex curve on the original hyperk\"ahler
manifold, which is impossible, because a subvariety
of a K\"ahler manifold cannot be homologous to zero.

In \ref{thm_Hermitian} we modify this argument
showing that it remains valid for any positive, exact
$(n-1, n-1)$-current in place of the current of
integration of $C$. 

A theorem of Sullivan \cite{_Sullivan_} implies that a complex manifold
admits no positive, exact $(n-1, n-1)$-currents
if and only if it is Hermitian symplectic, that is,
admits a taming symplectic structure (\ref{thm_currents}).

Streets and Tian have conjectured that Hermitian symplectic manifolds 
are K\"ahler, see \cite{_Streets_Tian_}. This conjecture was
treated at great length, and solved for some classes
of complex manifolds (for complex surfaces in \cite[Theorem 1.2]{_Li_Zhang_},
for twistor spaces in \cite{_MV:curves_twistors_}), but in full generality
it remains very difficult. In
\cite{_EFV:tamed_nilma_implies_Kahler_false_,_FV:erratum_} the
Streets--Tian conjecture was proven for complex
nilmanifolds.

The Streets--Tian conjecture makes sense for almost complex
manifolds as well; in this context, it is known as the 
``Donaldson's tamed to compatible conjecture''.
In \cite[Question 2]{_Donaldson:2-forms_}, Donaldson asked whether an almost
complex 4-manifold admitting a taming symplectic
form must also admit a compatible symplectic form.

Using \ref{thm_currents}, we reduce the proof of
K\"ahlerness of degenerate twistor deformations to
the Streets--Tian conjecture, and in Subsection 
\ref{_non-separa_Hermitian_sympl_Subsection_}, we prove
a special case of this conjecture, finishing the proof
of K\"ahlerness.

\hfill

\theorem
Let $\Teich$ be the Teichm\"uller space of holomorphically
symplectic complex structures, and $I\in \Teich$ 
a holomorphically symplectic complex structure
obtained as a limit $I= \lim_i I_i$, where $I_i\in \Teich$
are holomorphically
symplectic complex structures of K\"ahler type.
Assume that $(M, I)$ is Hermitian symplectic.
Then $(M,I)$ is also of K\"ahler type.

\hfill

\proof By \ref{_two_fibers_Herm_symple_bimero_Proposition_}, 
$(M,I)$ is of Fujiki class C. By \cite[Theorem 0.2]{Ch}, any Hermitian symplectic
Fujiki C manifold is of K\"ahler type.
\endproof


\section{Preliminaries}


\subsection{Hyperk\"ahler manifolds}

We start by recalling the basic facts of the theory of hyperk\"ahler manifolds.
For a detailed exposition of this theory see \cite{Hu}.

We will denote by $\bbH$ the algebra of quaternions and by $Sp(n)$
the group of quaternionic-linear transformations of $\bbH^n$.
Let $M$ be a compact $C^\infty$ mani\-fold. A Riemannian metric $g$ on $M$ is called {\bf hyperk\"ahler}
if the holo\-nomy group of its Levi--Civita connection $\nabla^g$ is contained in $Sp(n)$.
The elements of $Sp(n)$ commute with a two-dimensional family of complex structures
defined via multiplication by imaginary unit quaternions. Hence by the holo\-nomy principle
we get a family of $\nabla^g$-parallel complex structures on $M$ called the {\bf twistor family}.

In this paper we will always assume that hyperk\"ahler metrics are of {\bf maximal holonomy},
i.e. the holonomy group of $\nabla^g$ is isomorphic to $Sp(n)$. 
Given a complex structure $I$ on $M$, we will say that $I$ is of {\bf hyperk\"ahler type}
if there exists a hyperk\"ahler metric of maximal holonomy $g$ on $M$ such that $\nabla^g I = 0$.
If $I$ is of hyperk\"ahler type and $M$ is simply connected, the manifold $(M, I)$ is usually called irreducible
holomorphic symplectic (IHS).

Given a hyperk\"ahler metric $g$ on $M$ one can find $\nabla^g$-parallel complex structures $I$, $J$ and $K$
such that $K = IJ=-JI$. Denote by $\omega_I$, $\omega_J$ and $\omega_K$ the corresponding K\"ahler forms.
The 2-form $\Omega = \omega_J +\ii\omega_K$ is a holomorphic symplectic form on $(M, I)$.

Recall the following facts about the cohomology of a hyperk\"ahler mani\-fold $M$ of maximal holonomy
(see \cite{_Beauville_} and \cite{Hu} for the details).
The vector space $H^2(M,\bbQ)$ carries a non-degenerate quadratic form $q$
called the Beauville--Bogomolov--Fujiki (BBF) form. It has the following property:
there exists a positive constant $c_M\in \bbQ$
such that for all $a\in H^2(M,\bbQ)$ the {\bf Fujiki relation} holds:
\begin{eqnarray}\label{eqn_Fujiki}
\int_M a^{2n} = c_M q(a)^n.
\end{eqnarray}
We will normalize $q$ in such a way that it is primitive and integral on $H^2(M,\bbZ)$.
The BBF form has signature $(3,b_2(M)-3)$.

Consider the collection $\Comp(M)$ of complex structures of hyperk\"ahler type on $M$
and the group $\Diff^\circ(X)$ of diffeomorphisms of $M$ isotopic to the identity.
The quotient $\Teich(M) = \Comp(M)/\Diff^\circ(M)$ carries the structure of a non-Hausdorff complex
manifold and is called the {\bf Teichm\"uller space} of $M$.
The global Torelli theorem for hyperk\"ahler manifolds \cite{_Verbitsky:Torelli_}
describes the Teichm\"uller space in terms of the cohomology of $M$.
Define the {\bf period domain}:
$$
\DD = \{x\in \bbP H^2(M,\bbC)\st q(x) = 0,\, q(x,\bar{x})>0\}
$$
If $I$ is a complex structure of hyperk\"ahler type then the space
of $I$-holomor\-phic two-forms on $M$ is spanned by the symplectic form $\Omega$,
and the corresponding point $[\Omega]\in \DD$ is called the period of $I$.
This defines the {\bf period map} $\Per\colon \Teich(M)\to \DD$.
The global Torelli theorem \cite{_Verbitsky:Torelli_} claims that the restriction of $\Per$ to every
connected component of $\Teich$ is the Hausdorff reduction map, i.e. after
identifying the non-separated points of $\Teich$ it induces
an isomorphism of complex manifolds.

\subsection{Degenerate twistor deformations}

In this subsection, we 
briefly recall the definition of a degenerate twistor deformation associated
with a holomorphic Lagrangian fibration on a hyperk\"ahler manifold,
following \cite{BDV} (see also \cite{_V:degenerate_} and \cite{SV}).  
Holomorphic Lagrangian fibrations are ubiquitous in hyperk\"ahler
geometry, as shown by the following foundational theorem.

\hfill

\theorem (D. Matsushita)
Let $\pi\colon M \to B$ be a surjective holomorphic map
from a compact hyperk\"ahler manifold $M$ of maximal holonomy to a complex variety $B$ with $0<\dim B < \dim M$.
Then $\dim B = {1\over 2} \dim M$, and the fibers of $\pi$ are 
holomorphic Lagrangian.

\proof \cite{_Matsushita_}. \endproof 

\hfill

A map as in the above theorem is called a {\bf holomorphic Lagrangian fibration}.
The base $B$ is conjectured to be $\C P^n$, and this is the case in all
known examples; when $B$ is smooth, it is biholomorphic to $\C P^n$,
as shown in \cite{_Hwang:base_}.

\hfill

\definition\label{_C_symp_Definition_}
Let $M$ be a real $4n$-dimensional $C^\infty$-manifold.
A {\bf C-symplectic form} on $M$ is a smooth section $\Omega$
of the complexified bundle of 2-forms $\Lambda^2 M\otimes \bbC$
such that: $d\Omega = 0$, $\Omega^{n+1} = 0$ and the form
$\Omega^n\wedge\overline{\Omega}\mathstrut^n$ is pointwise non-vanishing.

\hfill

Let $\Omega$ be a C-symplectic form on $M$. For a point $x\in M$
let $T^{0,1}_x M \subset (T_x M)\otimes \bbC$ be the {\bf kernel}
of $\Omega$, i.e. the subspace consisting of complex tangent vectors $v\in (T_x M)\otimes \bbC$
such that the contraction $v\cnv\Omega$ vanishes.

One verifies (see \cite{BDV}, \cite{SV}) that the kernel of $\Omega$
has rank $2n$ at all points of $M$,
hence $T^{0,1}M \subset TM\otimes\bbC$ is a subbundle.
Defining $T^{1,0}M$ as the complex conjugate of $T^{0,1}M$,
one gets the decomposition $TM\otimes\bbC = T^{1,0}M\oplus T^{0,1}M$.
Using the closedness of the form $\Omega$, one checks that
the subbundle $T^{1,0}M$ is closed under the Lie bracket.
This defines an integrable almost complex structure on $M$,
and therefore by the Newlander--Nirenberg theorem
any manifold with a C-symplectic form admits an intrinsically defined
complex structure.

\hfill

\theorem\label{thm_csympl_lagr}
Let $M$ be a compact complex manifold with a holomorphic symplectic form $\Omega$,
and let $\pi\colon M \to B$ be a holomorphic Lagrangian fibration.
Consider a closed two-form 
$\eta\in \Lambda^2 B\otimes\bbC$ such that $\eta^{0,2}=0$. Then the two-form $\Omega_t:=\Omega+t \pi^*\eta$
is C-symplectic for any $t\in\bbC$. If we denote by $I_t$ the complex
structure on $M$ corresponding to $\Omega_t$, then for any $t\in \bbC$ the map $\pi$ is
holomorphic with respect to $I_t$ and the fixed complex structure on the base $B$.

\proof \cite{BDV}. \endproof

\hfill

The following property of the complex structures $I_t$ from the above theorem
appears implicitly in \cite{BDV} and \cite{SV}. Let us record it for the later use.

\hfill

\lemma\label{lem_invariance} In the setting of \ref{thm_csympl_lagr} denote by $F$ one of the fibres
of $\pi$ and let $x\in F$. Let $V = T_x F\otimes \bbC$ be the complexified tangent space
of $F$ at the point $x$. Then the subspace $V_t^{0,1} = V\cap T^{0,1}_{I_t} M$
does not depend on $t\in \bbC$. 

\hfill

\proof
Consider $v\in V$. Then $v\in V^{0,1}_t$ if and only if $v\cnv \Omega_t = 0$.
Note that $V = \ker(d\pi)_x$, therefore $v\cnv \pi^*\eta = 0$ and $v\cnv \Omega_t = v\cnv \Omega$.
Hence $V^{0,1}_t = V^{0,1}_0$, so $V^{0,1}_t$ is a fixed subspace of $V$
not depending on $t$.
\endproof

\hfill

Assume now that $M$ is a compact simply connected hyperk\"ahler manifold of maximal holonomy,
and $\pi\colon M \to B$ a holomorphic Lagrangian fibration over a projective base $B$
(the projectivity of $B$ follows from the other assumptions, see 
\cite[footnote on page 53]{_Amerik_Campana:families_}).
Let $\eta\in \Lambda^{1,1}B$ be the restriction of the Fubini--Study form under
some projective embedding of the base, $\Omega_t = \Omega+t\pi^*\eta$
the C-symplectic form and $I_t$ the corresponding complex structure, where $t\in \C$.
The manifolds $\MM_t = (M, I_t)$ are fibres of the {\bf the degenerate twistor deformation} of $M$,
which is a smooth holomorphic family $f\colon \MM\to \bbA^1$ over the affine line, see \cite{_V:degenerate_}.
This deformation has a number-theoretic interpretation (``the Shafarevich--Tate family''),
see \cite{_Markman:Lagrangian_,_AbashevaRogov_}.

The periods of the manifolds $\MM_t$
form an affine line in the period domain. More precisely, consider the subspace $W\subset H^2(M,\bbC)$
spanned by $\Omega$, $\bar{\Omega}$ and $\pi^*\eta$. Then the restriction of the BBF form $q|_W$ is
a degenerate positive-semidefinite quadratic form with one-dimensional kernel spanned by $\pi^*\eta$.
The corresponding quadric $Q_W$ is the union of two projective lines intersecting at the point $[\pi^*\eta]$.
The complement to that point in $Q_W$ is the disjoint union of two affine lines, one of them
passing through $[\Omega]$ and the other through $[\bar{\Omega}]$. The first of the two affine lines
is formed exactly by the periods of the manifolds $\MM_t$, $t\in\bbC$ of the degenerate twistor
family. The other affine line corresponds to the complex-conjugate degenerate twistor family
which one obtains replacing $\Omega$ by $\bar{\Omega}$.

\subsection{Weakly positive and strongly positive forms}

When we speak about ``weakly'' and ``strongly positive''
forms or currents, we use the French convention
about positivity: ``positive'' always  means $\geq 0$.
We follow \cite{_Demailly:analytic_,_Harvey_}.

\hfill

\definition
Let $(V,I)$ be a $2n$-dimensional real vector space equipped with a complex structure
operator $I$, $I^2=-\Id$.
A {\bf weakly positive $(p,p)$-form} on $V$ is a real $(p,p)$-form 
$\eta\in \Lambda^{p,p}_\bbR(V^*)$ which satisfies \[ 
\eta(x_1,Ix_1,x_2,Ix_2,\ldots, x_{p},Ix_p)\geq 0
\]
for all $x_1,\ldots, x_p \in V$. 

\hfill

\remark Clearly, the set of weakly positive $(p,p)$-forms 
is a convex cone. 

\hfill

\definition
The cone of {\bf strongly positive} polyvectors is the cone in $\Lambda^{p,p}_\bbR(V)$
generated by the monomials $x_1\wedge Ix_1\wedge x_2\wedge Ix_2\wedge \ldots\wedge x_p\wedge Ix_p$
with positive coefficients. Replacing $V$ with $V^*$, we also define
{\bf weakly positive $(p,p)$-polyvectors} and 
{\bf strongly positive $(p,p)$-forms.}

\hfill

\remark By definition, the cone of strongly positive polyvectors
is dual to the cone of weakly positive forms under the
natural pairing $\Lambda^{p,p}_\bbR(V^*)\times \Lambda^{p,p}_\bbR(V)\arrow \R$.

\hfill

\remark
Fix a volume element $\Vol\in \Lambda^{2n}(V^*)$ positive with respect to $I$.
For any $v\in \Lambda^{m}(V)$, consider the contraction
$v\cnv\Vol\in \Lambda^{2n-m}(V^*)$. For any monomial $\eta\in \Lambda^{m}(V)$,
$\eta\cnv\Vol$ is the complementary monomial.
Therefore, the contraction takes the strongly positive
cone in $\Lambda^{p,p}_\bbR(V)$ to the strongly positive cone
in $\Lambda^{n-p,n-p}_\bbR(V^*)$. Using this remark, it is not hard to prove the following basic claim.

\hfill

\proposition\label{prop_positive}
Consider the natural non-degenerate pairing
$$\Lambda^{p,p}(V^*)\times \Lambda^{n-p,n-p}(V^*)\to \R, \quad (\alpha, \beta)\mapsto \frac{\alpha \wedge \beta}{\Vol}.$$
Under this pairing the cone of strongly positive forms is dual to the cone of weakly positive forms.
Moreover, if $\alpha\in\Lambda^{p,p}(V^*)$ is weakly positive and $\beta\in\Lambda^{q,q}(V^*)$ is
strongly positive, then $\alpha\wdg\beta$ is weakly positive.

\proof \cite{_Demailly:analytic_}.
\endproof

\hfill

\lemma\label{lem_weaklypos}
Let $\eta$ be a $(p,0)$-form. Then
the form $(-1)^{\frac{p(p-1)}2}(\1)^p\eta\wedge\bar\eta$ is weakly positive.

\hfill

\proof
Clearly, a $(p,p)$-form $\alpha \in \Lambda^{p,p}(V^*)$ is weakly positive 
iff its restriction to any $p$-dimensional $I$-invariant subspace
$W\subset V$ is a non-negatively oriented volume form.
However, $\eta\restrict W$ is a complex linear volume form,
hence it can be written as $u z_1\wedge z_2\wedge \ldots \wedge z_p$,
where $z_i \in \Lambda^{1,0}(W^*)$ is a basis and $u\in\bbC$. 
Write $z_i:= x_i + \1 Ix_i$, where $x_i\in W^*$.
Then $z_i \wedge \bar z_i = -2\1 x_i \wedge I x_i$ which
gives 
\[ 
u z_1\wedge z_2\wedge ... \wedge z_p\wedge \bar u \bar z_1\wedge \bar z_2\wedge ... \wedge \bar z_p
= (-1)^{\frac{p(p-1)}2} (-2\1)^p |u|^2\prod_{i=1}^p x_i \wedge I x_i.
\]
The volume form $\prod_{i=1}^p x_i \wedge I x_i$ is positively oriented.
This implies that $(-1)^{\frac{p(p-1)}2}(\1)^p\eta\wedge\bar\eta$ is weakly positive.
\endproof

%

\subsection{Positive currents}

Let $M$ be a compact complex manifold of complex dimension $n$. Denote by $\EE^k(M)$
the space of complex $k$-forms of class $C^{\infty}$ on $M$. Fixing an arbitrary
Hermitian metric and an affine connection on $M$, one may
define the $C^d$-norm of $k$-forms, for arbitrary $d\ge 0$.
It is well known that this family of norms makes the space of smooth $k$-forms $\EE^k(M)$
a Fr\'echet topological vector space, the topology being independent
on the Hermitian metric and the connection. 
The space $\EE^k(M)$ has the Hodge decomposition
$$
\EE^k(M) = \bigoplus_{p+q=k} \EE^{p,q}(M), 
$$
where $\EE^{p,q}(M)$ are the closed subspaces of differential
forms of type $(p,q)$. We will denote by $\EE^k(M,\bbR)$ the space
of real $k$-forms with the induced topology.

Denote by $\EE'_k(M)$ the topological dual of $\EE^k(M)$, i.e.
the space of continuous linear forms on $\EE^k(M)$. We will always
work with the weak topology on $\EE'_k(M)$, i.e. the coarsest
topology for which all the evaluation maps $\mathrm{ev}_\alpha\colon \EE'_k(M)\to\bbC$,
$T\mapsto T(\alpha)$ for arbitrary $\alpha\in\EE^k(M)$ are continuous.
The topological vector space $\EE'_k(M)$ is called the space of
{\bf currents} of dimension $k$. The pairing
$\EE^{2n-k}(M)\times \EE^k(M)\to \bbC$ given by the wedge product and
integration over $M$ induces a natural continuous
embedding $\EE^{2n-k}(M)\hrarr \EE'_k(M)$. For this reason the
elements of $\EE'_k(M)$ are also called currents of {\bf degree $2n-k$}.

The Hodge decomposition of $k$-forms induces the Hodge decomposition
of currents. We denote by $\EE'_{p,q}(M)$ the dual of $\EE^{p,q}(M)$
and call it the space of currents of bidimension $(p,q)$ or of bidegree
$(n-p,n-q)$. The space of real currents will be denoted by $\EE'_k(M,\bbR)$
and $\EE'_{p,p}(M,\bbR) = \EE'_{p,p}(M)\cap \EE'_{2p}(M,\bbR)$.

The de Rham differential $d\colon \EE^k(M)\to \EE^{k+1}(M)$ is continuous,
and by the well known theory of elliptic complexes its image is closed in Fr\'echet topology.
It follows from the Banach open mapping theorem and the Hahn--Banach theorem
that the adjoint map $d\colon \EE'_{k+1}(M)\to \EE'_k(M)$ also has closed
image, i.e. the space of exact currents is closed in the weak topology.

Recall that a real current $T\in \EE'_{1,1}(M,\bbR)$ is called {\bf positive}
if for any $\eta\in\EE^{1,0}(M)$ we have $\ii T(\eta\wdg\bar{\eta})\ge 0$.
Let us denote the set of positive currents by $\EE^+_{1,1}(M)$. This set
is clearly a closed convex cone in $\EE'_{1,1}(M,\bbR)$. The positive currents
admit the following description, see e.g. \cite[\S 3]{HL}. Given a positive current
$T\in \EE^+_{1,1}(M)$ there exists a positive Borel measure $\|T\|$ called
the total variation of $T$ and a $\|T\|$-measurable section $\tau$ of the bundle
$\Lambda^{1,1}TM$, i.e. a $\|T\|$-measurable bivector field of Hodge type $(1,1)$,
such that $T = \tau \|T\|$. The latter equality means that the action of $T$
on $\alpha\in\EE^{1,1}(M)$ is given by
$$
T(\alpha) = \int_M \langle\tau,\alpha\rangle d\|T\|,
$$
where $\langle\tau,\alpha\rangle$ denotes the pairing between bivectors and two-forms.
The bivector field $\tau$ is positive in the following sense:
for any $\eta\in\EE^{1,0}(M)$ the function $\langle \tau, \ii\eta\wdg\bar{\eta}\rangle$
is non-negative $\|T\|$-almost everywhere on $M$. 

\hfill

\proposition\label{prop_tangency}
Let $\pi\colon M \to X$ be a proper holomorphic map,
and $T\in \EE^+_{1,1}(M)$
a positive current such that $T(\pi^*\alpha)=0$
for any positive $(1,1)$-form $\alpha$ on $X$.
Then in the above decomposition $T = \tau\|T\|$ the
bivector field $\tau$ is almost everywhere tangent
to the fibers of $\pi$, i.e. $\tau \in \ker{d\pi}$
outside a subset of $\|T\|$-measure zero.

\hfill

\proof 
The differential $d\pi:\; T_x M \arrow T_{\pi(x)}X$
can be interpreted as a map from $TM$ to $\pi^* (TX)$. Then 
$d\pi(\tau)$ is a measurable section of the bundle $\pi^*\Lambda^{1,1}TX$.
By our assumptions, for any $\eta\in\EE^{1,0}(X)$ we have
$$
\ii T(\pi^*\eta\wdg\pi^*\bar{\eta}) = \int_M \langle d\pi(\tau), \ii\eta\wdg\bar{\eta}\rangle d\|T\| = 0,
$$
and since $d\pi(\tau)$ is positive, the function $\langle d\pi(\tau), \ii\eta\wdg\bar{\eta}\rangle$
is non-negative almost everywhere on $M$, therefore it must vanish almost everywhere.
Since this is true for arbitrary $\eta$ as above, we conclude that $d\pi(\tau)$
vanishes almost everywhere, which concludes the proof.
\endproof

\section{Hermitian symplectic structures}

In the sequel, we use the following version of the Hahn--Banach theorem,
sometimes called the ``Hahn--Banach separation theorem''.

\hfill

\theorem\label{_H_B_Separation_Theorem_}
(Hahn--Banach theorem, \cite{_Bourbaki:TVP_}, or \cite[Chapter II, Theorem 3.1]{_Schaefer_}) \\
Let $V$ be a topological vector space, $A\subset V$ an open convex subset
of $V$, and $W$ a subspace of $V$ satisfying $W\cap A=
\emptyset$. Then there exists a continuous linear functional
$\theta$ on $V$, such that
$\theta\restrict A>0$ and $\theta\restrict W = 0$.

\hfill

Recall the following definitions.

\hfill

\definition
Let $(M, I)$ be a complex manifold and $\omega\in \Lambda^2 M$ a real symplectic form.
The complex structure
$I$ is {\bf tamed} by $\omega$ if for any non-zero tangent vector $v$ at an arbitrary
point of $M$ we have $\omega(v,Iv)> 0$. In this case we call $\omega$ a
{\bf Hermitian symplectic structure.} Equivalently, a symplectic
form $\omega$ is Hermitian symplectic if its $(1,1)$-part $\omega^{1,1}$
is a Hermitian form.

\hfill

The following characterization 
of compact complex manifolds admitting Hermitian symplectic
structures is originally due to Sullivan \cite{_Sullivan_}, see \cite{Ch}, \cite{_DP:hermitian_} and \cite{HL}
for related results and a broader context.

\hfill

\theorem\label{thm_currents}
Let $M$ be a compact complex manifold. The manifold $M$ 
does not admit a Hermitian
symplectic structure if and only if 
if admits a non-zero positive exact 
$(n-1,n-1)$-current.

\hfill

\proof 
We sketch the well known proof for the sake of completeness
and because we refer to these arguments elsewhere.
We use the Hahn--Banach separation theorem 
(\ref{_H_B_Separation_Theorem_}), applied to $V$, $A$ and $W$
defined below.

Let $V=\Lambda^2(M,\R)$ be the Fr\'echet space of smooth 
two-forms equipped with the $C^\infty$-topology, $A$ the cone 
of two-forms with strictly positive $(1,1)$-part,
and $W\subset V$ the space of closed two-forms.
The intersection $A \cap W$ is the set of all 
Hermitian symplectic forms. If this set is empty, we 
use \ref{_H_B_Separation_Theorem_} to find a 
non-zero functional $\theta\in V^*$ which is positive on $A$
(this means that $\theta$ is a positive 
$(n-1,n-1)$-current) and vanishing on $W$. 
The latter condition means that $\theta$ is exact.
Indeed, $\langle \theta, du\rangle=0$ for all $u\in \Lambda^1 M$
implies that $\theta$ is closed. It is exact, because otherwise
it is non-zero on $W$ by Poincar\'e duality (the de Rham cohomology
of currents is equal to the de Rham cohomology of differential
forms, \cite{_Demailly:analytic_}).
\endproof

\hfill

\theorem\label{thm_Hermitian}
Let $\pi\colon M \to B$ be a holomorphic Lagrangian
fibration on a compact simply connected hyperk\"ahler manifold
of maximal holonomy, and $\MM_t = (M, I_t)$, $t\in \C$ the corresponding 
degenerate twistor deformation.
Then $\MM_t$ admits an Hermitian symplectic structure
for all $t\in \C$.

\hfill

\begin{proof}
We use \ref{thm_currents}. Let $T\in \EE^+_{1,1}(\MM_t)$ be an exact
positive current of type $(1,1)$ on $\MM_t$. We need to prove that
$T = 0$.

As explained above, $T$ admits a representation $T = \tau\mu$, where
$\mu = \|T\|$ is a positive Borel measure and $\tau$ is a $\mu$-measurable
bivector field on $\MM_t$. Let us check that the conditions of \ref{prop_tangency}
are satisfied. Fix a K\"ahler form $\eta\in\Lambda^{1,1}B$. For any positive
form $\alpha\in \Lambda^{1,1}B$ we have $0\le \alpha\le C\eta$ for some constant $C>0$.
Then $0\le T(\pi^*\alpha) \le T(\pi^*\eta)$, but $\pi^*\eta$ is closed, and
since $T$ is exact $T(\pi^*\eta) = 0$. Therefore $T(\pi^*\alpha) = 0$ and
we conclude by \ref{prop_tangency} that the bivector field $\tau$ is tangent
to the fibers of $\pi$ (after, possibly, changing $\tau$ on a subset of $\mu$-measure zero).

Fix a point $x\in M$, let $V = T_x M\otimes \bbC$ and let $W=\ker(d\pi)_x\subset V$ be the
complexified tangent space to the fiber of $\pi$ through $x$. The bivector $\tau(x)$
is a linear combination with positive coefficients of bivectors of the form
$-\ii u\wdg \bar{u}$ for $u\in V^{1,0}_t$. But since $\tau$ is tangent to the
fibers of $\pi$, we have $u\in W^{1,0}_t$ for all bivectors appearing in that linear
combination. Now, we know by \ref{lem_invariance} that the subspace $W^{1,0}_t \subset W$
does not depend on $t$. We conclude that the bivector field $\tau$ is of type $(1,1)$
and positive for all complex structures $I_t$ in the degenerate twistor family.
Therefore $T$ is an exact positive current of bidimension $(1,1)$ on $\MM_0 = M$.
The latter being a K\"ahler manifold, the current $T$ must vanish by \ref{thm_currents}.
\end{proof}


\section{$d$-Commendable manifolds}


\definition Let $M$ be a compact complex manifold, $h\in \Lambda^{1,1}M$ a Hermitian
form and $0 < d < \dim_\bbC(M)$.
We say that $M$ is {\bf $d$-commendable} if there exists $\eta \in \Lambda^{2d} M$
such that $d\eta = 0$ and the form 
$\eta^{d,d}- \epsilon h^d$ is weakly positive for some $\epsilon>0$.\footnote%
{In the sequel, we will denote this condition as $\eta^{d,d}\geq \epsilon h^d$ .}
Note that 1-commendable is the same as Hermitian symplectic.

\hfill

\proposition\label{prop_commendable}
Let $M$ be a complex manifold with a Hermitian symplectic structure $\omega$.
Then $M$ is $d$-commendable, for all $d> 0$.

\hfill

\begin{proof}
For $0 < d < \dim_\bbC(M)$ let $\eta =\omega^d$. Then $\eta$ closed and
$\eta^{d,d}$ is a positive linear combination of the forms $(\omega^{2,0} \wdg \overline{\omega^{2,0}})^k \wdg (\omega^{1,1})^{d-2k}$,
for $k=0,\ldots, [d/2]$. The form $(\omega^{2,0})^k\wdg (\overline{\omega^{2,0}})^k$ is weakly positive
by \ref{lem_weaklypos},
and $(\omega^{1,1})^{d-2k}$ is strongly positive for $2k<d$. By \ref{prop_positive}
the forms $(\omega^{2,0} \wdg \overline{\omega^{2,0}})^k \wdg (\omega^{1,1})^{d-2k}$ are weakly positive. Moreover, $(\omega^{1,1})^d \ge Ch^d$ for
some constant $C>0$, because $\omega^{1,1}$ is a Hermitian form. Therefore $\eta^{d,d}\ge C h^d$.
\end{proof}

\subsection{Barlet spaces on commendable manifolds}

The purpose of introducing the notion of $d$-commendability is that it
can serve as a replacement of the condition of being K\"ahler in some situations.
More precisely, we will observe that the Barlet spaces of $d$-dimensional cycles
of $d$-commendable manifolds behave in some ways similar to the spaces of cycles
of K\"ahler manifolds. Let us start with the following claim.

\hfill

\proposition \label{_Bishop_Herm_Sympl_Proposition_}
Let $M$ be a $d$-commendable compact complex manifold and $\CC_d$ a connected
component of the Barlet space of dimension $d$ cycles on $M$. Then $\CC_d$ is compact.

\hfill

\begin{proof}
Let $h\in \Lambda^{1,1}M$ be a Hermitian metric on $M$ and $\eta$ a closed $2d$-form with $\eta^{d,d}\ge C h^d$.
By the compactness criterion of Barlet, based on Bishop's theorem
(\cite{_Magnusson:cycle_}, \cite[Theorem 4.2.69]{BM}), we need to show
that the volume of all cycles in $\CC_d$ with respect to $h$ is bounded.
Assume that $[Z]\in \CC_d$ is a cycle. Then
$$
\int_Z h^d \le (1/C)\int_Z \eta^{d,d} = (1/C)\int_Z \eta.
$$
Since $d\eta = 0$, the latter integral depends only on the cohomology class of $[Z]$ which is
constant for all cycles in the connected component $\CC_d$. The claim follows. 
\end{proof}

\hfill

In what follows we will need a relative version of the above statement.
We will consider the following setting. Let $\pi\colon \MM\to B$ be a proper
smooth holomorphic map between complex manifolds. Let $\CC_d(\MM/B)$ be a connected
component of the Barlet space of $d$-dimensional cycles contained in the fibers
of $\pi$, see \cite[section 4.8.2]{BM}. Then the map $p\colon \CC_d(\MM/B)\to B$ that sends a cycle to
the point $t\in B$ over which this cycle is supported is holomorphic (loc. cit.)
The following proposition is a generalization of \cite[Corollary 4.8.]{BM}.

\hfill

\proposition\label{prop_proper}
In the above setting assume that there exists a differential form $\eta\in \Lambda^{2d}(\MM)$
such that for any $t\in B$ the fiber $\MM_t$ and the restriction $\eta_t = \eta|_{\MM_t}$
satisfy the conditions of $d$-commendability. Then the support morphism $p\colon \CC_d(\MM/B)\to B$
is proper, in particular its image is an analytic subvariety of $B$.

\hfill

\begin{proof}
To prove properness we restrict to a compact subset $K\subset B$.
The cycles in the preimage $p^{-1}(K)$ are supported in the compact
subset $\pi^{-1}(K)$ of $\MM$. So by the compactness criterion
\cite[Theorem 4.2.69]{BM} we need to prove that the volume of the cycles
$[Z]\in \CC_d(\MM/B)$ in $p^{-1}(K)$ is bounded above. Let us fix a Hermitian metric
$h$ on $\MM$. 


We can assume that $K$ is sufficiently small and find a constant $C>0$ such that for any $t\in K$ we have
$\eta_t^{d,d}\ge C h_t^d$, where the subscript $t$ denotes the restriction
to the fiber $\MM_t$. Now, given a cycle $[Z]\in p^{-1}(K)$ let $t\in K$ be
the point over which $[Z]$ is supported. We estimate
the volume of $[Z]$:
$$
\int_Z h^d = \int_{Z} h^d_t \le (1/C)\int_Z \eta_t^{d,d} = (1/C)\int_Z \eta_t.
$$
The homology class of $[Z]$ is constant, and the restriction of $\eta$
to any fiber $\MM_t$ is closed, so the integral $\int_Z \eta_t$ depends
only on the point $t\in K$ over which $[Z]$ is supported; let us
denote the value of the integral by $\phi(t)$. Since
$\eta$ is a smooth differential form on $\MM$, $\phi$ is a continuous function. Since $K$ is compact,
$\phi$ is bounded from above on $K$. This gives the required upper bound for the volume of
the cycles $[Z]\in p^{-1}(K)$.
\end{proof}

\subsection{Non-separated points in the Teichm\"uller space}
\label{_non-separa_Hermitian_sympl_Subsection_}

Let $M$ be a simply connected compact hyperk\"ahler manifold of maximal holonomy
and $\pi\colon M\to B$ a holomorphic Lagrangian fibration. Let $\MM\to\bbA^1$
be the corresponding degenerate twistor family. Denote by $U\subset \bbA^1$ the
open subset of points of the base such that for $t\in U$ the fiber $\MM_t$ is K\"ahler.

\hfill

The argument used to prove
\ref{_two_fibers_Herm_symple_bimero_Proposition_} below
is a version of the one used in \cite{Hu}
to establish the bimeromorphic correspondence
between non-separated fibers of a twistor deformation.
In  \cite{Hu}, both fibers are assumed to be K\"ahler.
We prove the same statement when the fibers are Hermitian symplectic.

\hfill

By the global
Torelli theorem (\cite{_Verbitsky:Torelli_}) 
the hyperk\"ahler manifold $\MM_t$ admits a hyperk\"ahler deformation $X$
such that the period of $X$ corresponds to the point $t_0$ and $X$ lies in the same
connected component of the Teichm\"uller space as $\MM_t$
(\cite[Section 1.2]{_V:degenerate_}). Let $\XX\to \Delta$ denote the
restriction of the universal deformation of $X$ to $\Delta$.

\hfill

\proposition\label{_two_fibers_Herm_symple_bimero_Proposition_}
Consider a family $\MM\to\bbA^1$ of holomorphically
symplectic manifolds, and let $U\subset \bbA^1$ be the set
of all points corresponding to K\"ahler fibers. 
Assume the the fiber over a limit point $t_0\in \dd U$ is
Hermitian symplectic.
Then there exists a hyperk\"ahler deformation $X$ of $M$
that is bimeromorphic to $\MM_{t_0}$

\hfill

\begin{proof}
Let $\Delta\subset \bbA^1$ be a small disc around $t_0$. Note that $\Delta\cap U\neq\emptyset$,
so we may find a K\"ahler fiber $\MM_t$ with
$t\in\Delta$. 

Let $\omega$ be a Hermitian symplectic form on
$\MM_{t_0}$, and $\omega'$ a K\"ahler form on $X$.

Note that $\omega+\omega'$ is a Hermitian symplectic form on $X\times \MM_{t_0}$.
By continuity, after possibly shrinking $\Delta$, we may assume that the restriction
of $\omega+\omega'$ to any fiber $\XX_t\times \MM_t$, $t\in \Delta$ makes this fiber 
Hermitian symplectic.

The idea now is to take the graph $\Gamma_t$ of the bimeromorphism $\XX_t\dashrightarrow \MM_t$
and use \ref{prop_proper} to show that there is a  subsequence
$\{\Gamma_{t_i}\}$ converging to a subvariety $\Gamma_{t_0}$ in $\XX_{t_0}\times \MM_{t_0}$.
Then we use a topological argument to show that one of the irreducible components of
$\Gamma_{t_0}$ is a graph of a bimeromorphism $\XX_{t_0}\dashrightarrow \MM_{t_0}$.

Let us consider the open subset $\Delta\cap U$. It is non-empty because $t_0$ is a limit
point of $U$. The fibers of the family $\MM$ over this open subset are K\"ahler by the definition
of $U$. The points of the Teichm\"uller space corresponding to $\MM_t$ and $\XX_t$,
$t\in \Delta\cap U$ are either equal or non-separated by construction, hence by the theorem of Huybrechts \cite[Theorem 4.3]{Hu}
the manifolds $\MM_t$ and $\XX_t$ are bimeromorphic for any $t\in \Delta\cap U$.

Consider the complex manifold $\YY = \XX\times_{\Delta} \MM$. For $t\in \Delta\cap U$ the
graph of a bimeromorphic map between $\XX_t$ and $\MM_t$ is a subvariety $\Gamma_t \subset \YY$.
Let $\CC_{\dim_\bbC X}(\YY/\Delta)$ be an irreducible component of the relative Barlet
space that contains a point corresponding to one of these graphs. By \ref{prop_proper}
the image of $\CC_{\dim_\bbC X}(\YY/\Delta)$ in $\Delta$ is a closed subvariety,
hence it is either a countable set of points, or the whole disc $\Delta$.
Since the subset $\Delta\cap U$ is uncountable, one of the connected components
of the Barlet space must map surjectively onto $\Delta$.

We conclude that there exists a sequence of cycles
 $[\Gamma_{t_i}]$ in the fibers $\YY_{t_i}$,
where $t_i\in \Delta\cap U$, $t_i\to t_0$ for $i\to\infty$, 
giving the above bimeromorphic maps, contained
in the same irreducible component of the Barlet space and converging to a limit cycle $[\Gamma_{t_0}]$
in the fiber $X\times \MM_{t_0}$ over $t_0$. Then the homology classes of $[\Gamma_{t_i}]$
do not depend on $i$. Denote by $\Omega'$ a holomorphic symplectic form on $X$.
The bimeromorphic maps induced by $\Gamma_{t_i}$ map a holomorphic symplectic form on $\MM_{t_i}$
to a holomorphic symplectic form on $\XX_{t_i}$. Therefore by continuity $[\Gamma_{t_0}]$
seen as a correspondence acting on cohomo\-logy maps the cohomology class $[\Omega_{t_0}]$
on $\MM_{t_0}$ to the cohomology class $[\Omega']$ on $X$. We can therefore apply
\ref{prop_bimer} below to conclude that the cycle $[\Gamma_{t_0}]$ contains exactly one
irreducible component that is the graph of a bimeromorphic map between $\MM_{t_0}$ and $X$.
\end{proof}

\hfill

The proof of the following proposition is a version of \cite[proof of Theorem 4.3]{Hu}.
We use it to establish the bimeromorphic equivalence of two holomorphically symplectic
manifolds, the only difference from loc. cit. being that we do not assume one of the
manifolds to be K\"ahler.

We consider the following setting. Assume that $M$ is a compact complex manifold with
a holomorphic symplectic form $\Omega$ and $X$ is a compact hyperk\"ahler manifold of maximal holonomy
with a holomorphic symplectic form $\Omega'$. Assume that $\dim_\bbC(M) = \dim_\bbC(X) = 2n$
and that $[Z]$ is a $2n$-cycle on $X\times M$. Assume that $[Z]$ is of degree one both over $X$
and over $M$. Assume moreover that $[Z]_*\colon H^2(M,\bbC)\to H^2(X,\bbC)$ is
an isomorphism that sends $[\Omega]$ to $[\Omega']$. We can decompose the
cycle into a finite sum $[Z] = \sum a_i [Z_i]$, where $a_i>0$ and $Z_i\subset X\times M$
is an irreducible subvariety.

\hfill

\proposition\label{prop_bimer} In the above decomposition of $[Z]$ there is
exactly one irreducible component dominating both $M$ and $X$ and this component
is the graph of a bimeromorphic isomorphism between $M$ and $X$.

\hfill

\begin{proof} The degree of $[Z] = \sum a_i [Z_i]$ over $X$ is equal to $\sum a_i \mathrm{deg}_X(Z_i)$,
and since all $a_i$ are positive, there exists $i_0$ such that $\mathrm{deg}_X(Z_{i_0}) = 1$, $a_{i_0} = 1$
and $\mathrm{deg}_X(Z_i) = 0$ for $i\neq i_0$. Analogously, there exists $i_1$, such that 
$\mathrm{deg}_M(Z_{i_1}) = 1$, $a_{i_1}=1$ and $\mathrm{deg}_M(Z_i) = 0$ for $i\neq i_1$.
We need to prove that $i_0 = i_1$: in this case the subvariety $Z_{i_0} = Z_{i_1}$ is the graph
of a bimeromorphic map between $X$ and $M$.

Denote by $p\colon X\times M\to X$ and $q\colon X\times M\to M$ the two projections.
Let us consider the action of the correspondence given by $[Z]$ on the cohomology class
of the symplectic form $\Omega$. Recall that by our assumptions $[Z]_*[\Omega] = [\Omega']$.
For a connected component $Z_i$ of the cycle $Z$ denote by $\tilde{Z}_i$ its resolution of
singularities and by $\tilde{p}_i\colon \tilde{Z}_i \to X$ and $\tilde{q}_i\colon \tilde{Z}_i \to M$
the induced maps. Then $[Z_i]_*[\Omega]$  is the cohomology clases of the current $\tilde{p}_{i*}(\tilde{q}_i^*\Omega)$.
For $i\neq i_0$ the map $\tilde{p}_i$ is not dominant, and since $\tilde{q}_i^*\Omega$ is a holomorphic
two-form, the current $\tilde{p}_{i*}(\tilde{q}_i^*\Omega)$ vanishes.
Indeed, a closed (2,0)-current is holomorphic, hence smooth by elliptic
regularity, but the support of $\tilde{p}_{i*}(\tilde{q}_i^*\Omega)$ 
is contained in a proper subvariety.
 Now consider the case $i=i_0$.
Then $\tilde{p}_{i_0*}(\tilde{q}_{i_0}^*\Omega)$ is a holomorphic two-form on $X$, therefore it
must be proportional to $\Omega'$, because $X$ is hyperk\"ahler of maximal holonomy.
It follows that $\tilde{q}_{i_0}^*\Omega$ is generically a symplectic form, so $\tilde{q}_{i_0}$
must be dominant, hence $i_0 = i_1$.
\end{proof}

\subsection{Degenerate twistor deformations are K\"ahler}

We finish the paper by proving our main result.

\hfill

\theorem\label{thm_main} 
Let $M$ be a simply connected compact
hyperk\"ahler manifold of maximal holonomy and $\pi\colon M\to B$
a holomorphic Lagrangian fibration. Let $\MM\to\bbA^1$ be the corresponding
degenerate twistor family. Then all fibers $\MM_t$, $t\in\bbA^1$ are
K\"ahler.

\hfill

\begin{proof}
Denote by $U\subset \bbA^1$ the open subset of points of the base
such that for $t\in U$ the fiber $\MM_t$ is K\"ahler. Assume that $U\neq \bbA^1$
and let $x\in \dd U$. Then by \ref{_two_fibers_Herm_symple_bimero_Proposition_}, the fiber $\MM_t$
is of Fujiki class C. Moreover, by \ref{thm_Hermitian} the manifold
$\MM_t$ admits a Hermitian symplectic form $\omega$. Then $\omega^{1,1}$ is a Hermitian
metric and $dd^c(\omega^{1,1}) = 0$, therefore $\MM_t$ is K\"ahler by \cite[Theorem 0.2]{Ch}.%
\footnote{We thank the participants of 
the \href{https://mathoverflow.net/questions/471007}{mathoverflow discussion}
for this reference}
It follows that $U = \bbA^1$.
\end{proof}

\hfill

{\bf Acknowledgements:} 
We are grateful to Ekaterina Amerik, Alexey Golota, Evgeny Shinder,
Olivier Martin and Eyal Markman for interesting discussions of the subject.

\hfill

{

}

\noindent {\sc Andrey Soldatenkov\\
{\tt aosoldatenkov@gmail.com}\\
{\sc Universidade Estadual de Campinas
Departamento de Matem\'atica - IMECC
Rua S\'ergio Buarque de Holanda, 651
13083-859, Campinas - SP, Brasil}}

\hfill

\noindent {\sc Misha Verbitsky\\
{\sc Instituto Nacional de Matem\'atica Pura e
              Aplicada (IMPA) \\ Estrada Dona Castorina, 110\\
Jardim Bot\^anico, CEP 22460-320\\
Rio de Janeiro, RJ - Brasil\\
\tt  verbit@impa.br }
}

\end{document}